\documentclass[10pt,rotate]{amsart}
\usepackage[all]{xypic}
\usepackage{latexsym}
\usepackage{amssymb}
\usepackage{amsfonts}
\usepackage{amscd}
\usepackage{amsmath,amsthm}

\newtheorem{lemma}{Lemma}[section]

\newtheorem{rem}[lemma]{Remark}

\newtheorem{lem}[lemma]{Lemma}

\newtheorem{prop}[lemma]{Proposition}
\newtheorem{thm}[lemma]{Theorem}
\newtheorem{cor}[lemma]{Corollary}

{
}

\theoremstyle{definition}

\newcommand{\Ups}{\Upsilon}
\newcommand{\Upsol}{\bar{\Upsilon}}

\theoremstyle{remark}

\numberwithin{equation}{section}

\newenvironment{pf}{\noindent{\bf Proof.}}{\hfill $\square$\medskip}

\def\AA{{\mathbb A}}
\def\CC{{\mathbb C}} 
 
\def\FF{{\mathbb F}} 
\def\GG{{\mathbb G}}

\def\NN{{\mathbb N}}

\def\PP{{\mathbb P}}

\def\ZZ{{\mathbb Z}}

\def\phiol{{\bar \phi}}

\def\psiol{{\bar \psi}}

\def\Aol{{\bar A}}

\def\Mol{{\bar M}}
\def\cMol{{\bar{\cal M}}}

\def\0ol{{\bar 0}}
\def\1ol{{\bar 1}}
\def\2ol{{\bar 2}}
\def\ol2{{\bar 2}}
\def\3ol{{\bar 3}}
\def\4ol{{\bar 4}}
\def\5ol{{\bar 5}}
\def\6ol{{\bar 6}}
\def\7ol{{\bar 7}}
\def\8ol{{\bar 8}}
\def\9ol{{\bar 9}}

\def\bold0{{\bf 0}}
\def\bold1{{\bf 1}}
\def\bold2{{\bf 2}} 
\def\bold3{{\bf  3}}
\def\bold4{{\bf 4}}
\def\bold5{{\bf 5}}
\def\bold6{{\bf 6}}
\def\bold7{{\bf 7}}
\def\bold8{{\bf 8}}
\def\bold9{{\bf 9}}

\def\P2Skly{\PP^2_{Skly}}

\def\coker{\operatorname {coker}}

\def\Ext{\operatorname {Ext}}

\def\gr{\operatorname {gr}}

\def\Hom{\operatorname {Hom}}

\def\ker{\operatorname {ker}}

\def\th{\operatorname {th}}    % for writing n^{th}
\def\Tor{\operatorname {Tor}}

\def\dim{\operatorname{dim}}

\def\Ext{\operatorname{Ext}}

\def\sgr{{\sf gr}}

\def\Gr{{\sf Gr}}

\def\Hom{\operatorname{Hom}}

\def\id{\operatorname{id}}

\def\Pic{\operatorname{Pic}}
\def\Proj{\operatorname{Proj}}

\def\Spec{\operatorname{Spec}}

\def\ul1{\operatorname{\underline{1}}}

\def\G{\mathop{\underline{\underline{\it \Gamma}}}\nolimits}

\def\l{\leftarrow}

\def\d{\downarrow}

\def\a{\alpha}
\def\b{\beta}

\def\d{\delta}

\def\l{\lambda}

\def\D{\Delta}
\def\G{\Gamma}
\def\L{\Lambda}

\def\fa{{\mathfrak a}}
\def\fb{{\mathfrak b}}

\def\fm{{\mathfrak m}}
\def\fn{{\mathfrak n}}

\def\sP{{\sf P}}

\def\sT{{\sf T}}

\def\wtM{{\widetilde{ M}}}

\def\wtS{{\widetilde{ S}}}

\def\cal{\mathcal}

\def\cF{{\cal F}}

\def\cM{{\cal M}}
\def\cN{{\cal N}}
\def\cO{{\cal O}}

\def\cU{{\cal U}}

\def\cX{{\cal X}}

\def\cZ{{\cal Z}}

\def\coh{{\sf coh}}

\def\Qcoh{{\sf Qcoh}}

\def\dirlim{\mathop{\vtop{\baselineskip -100pt\lineskip -1pt\lineskiplimit 0pt
\setbox0\hbox{lim}\copy0\hbox to \wd0{\rightarrowfill}}}\limits}
\def\invlim{\mathop{\vtop{\baselineskip -100pt\lineskip -1pt\lineskiplimit 0pt
\setbox0\hbox{lim}\copy0\hbox to \wd0{\leftarrowfill}}}\limits}

\def\I11{{1 \kern -0.8pt \! \mbox{l}}}
\def\mumu{{\mu\kern-4.2pt\mu}}
\def\bfmu{{\mu\kern-4.2pt\mu}}
\def\2slash{\backslash \! \backslash}

\def\boxtimes{\setbox0\hbox{$\Box$}\copy0\kern-\wd0\hbox{$\times$}}

\begin{document}
  
\title[Pic and $K_0$ for toric DM stacks]{Computation of the Grothendieck and Picard groups of a toric DM stack $\cX$
by using a homogeneous coordinate ring for $\cX$}
\author{S. Paul Smith} 
\address{Department of Mathematics, Box 354350, Univ.
Washington, Seattle, WA 98195, USA}
\email{smith@math.washington.edu}

\thanks{The author was supported by NSF grants DMS-0245724 and DMS-0602347 }

\begin{abstract}
We compute the Grothendieck and Picard groups of a smooth toric DM stack by using
a suitable category of graded modules over a polynomial ring. The polynomial ring with a suitable grading and suitable irrelevant ideal functions as a homogeneous coordinate ring for the stack.

\end{abstract}
 
\maketitle 

\pagenumbering{arabic}

\section{Introduction}

\subsection{}
Toric Deligne-Mumford stacks  were introduced in the  seminal paper of Borisov, Chen, and Smith \cite{BCS}.
A toric DM stack is a particular kind of quotient stack $\cX:=[(V-Z)/G]$ for an abelian affine algebraic group 
$G$ acting on a scheme $V-Z$ that is the complement to a suitable union of $G$-stable subspaces $Z$
in a rational representation $V$ of $G$. 

The appendix to Vistoli's  paper \cite[Example 7.21]{V} shows that $\Qcoh \cX$ is equivalent to the category of 
$G$-equivariant sheaves on $V-Z$. The analysis of toric DM stacks in \cite{BCS}, and in subsequent papers such as the framework developed by Fantechi, Mann, and Naroni \cite{FMN}, is carried out in the context of  $G$-equivariant sheaves on  $V-Z$. In particular, the computation of the 
Grothendieck group $K_0(\cX)$ by Borisov and Horja \cite[Thm. 4.10]{BH} is couched in the 
language of $G$-equivariant sheaves.

In this paper we offer a different approach that is arguably more elementary.

 \subsection{}
 Let $\G$ be the rational character group of $G$. It is a finitely generated abelian group, and may have torsion.
It is well known that 
\begin{enumerate}
  \item 
  the category of $G$-equivariant $\cO_V$-modules is equivalent
 to $\Gr(\cO(V), \G)$, the category of $\G$-graded $\cO(V)$-modules, and
  \item 
  the category of $G$-equivariant $\cO_{V-Z}$-modules is equivalent to the quotient of 
 $\Gr(\cO(V), \G)$ by the localizing subcategory generated by the coherent  $G$-equivariant $\cO_{V}$-modules that are supported on $Z$. 
\end{enumerate} 
It follows that there is an equivalence of categories
\begin{equation}
\label{eq.equiv}
\Qcoh \biggl[{{V-Z}\over{G}}\biggr] \equiv {{\Gr(\cO(V),\G)}\over{\sT}}
\end{equation}
where $\sT$ is  the localizing subcategory described in (2).   
This  is an equivalence of monoidal categories: the tensor product of
$\cO_{\cX}$-modules corresponds to the tensor product of graded $\cO(V)$-modules where one
uses the tensor product grading, $\deg(M_\a \otimes N_\b)=\a+\b$.

The graded ring $(\cO(V), \G)$ is analogous to Cox's homogeneous coordinate ring of a toric variety.

 \subsection{}
In this paper we use the quotient of the graded module category in (\ref{eq.equiv}) 
to study a toric DM stack. 
We prove two results. We compute the Grothendieck group $K_0(\cX):=K_0(\coh \cX)$ of 
the category of coherent sheaves on $\cX$ 
in Theorem \ref{thm.K0} and the Picard group $\Pic (\cX)$ in Theorem \ref{thm.pic}. 
The result on the Picard group is proved in somewhat greater generality.

\subsection{}
In section \ref{ssect.hcrs} we introduce the notion of a homogeneous coordinate ring for certain stacks. 
It consists of a triple $(A,\D,\fa)$  consisting of a commutative ring graded by an arbitrary finitely generated
abelian group $\D$ and an irrelevant graded ideal $\fa$. The idea is to mimic Cox's homogeneous coordinate ring for toric varieties. The homogeneous coordinate ring for the stacks we are interested in is particularly accessible and allows one to avoid certain local-to-global arguments and to avoid some of the more technical aspects of stacks.

We extend the definition of a connected graded ring to this setting, i.e., when the grading group is arbitrary,
and use that to prove more widely applicable versions of some standard results for connected graded rings (modeled on those for commutative local rings). For example, we prove a version of Nakayama's lemma, and as a consequence show that every finitely graded 
projective module is a direct sum of free modules in a unique way (Lemma \ref{lem.projs}.

This allows us to show that for certain stacks $\cX$ every coherent $\cO_\cX$-module has a finite resolution in
which each term is a  finite direct sum of invertible $\cO_\cX$-modules (Prop. \ref{prop.res}).

 \subsection{}
One of the main motivations for the introduction and study of toric DM stacks is that they are are a source of 
nice ambient spaces for the orbifolds that arise in string theory. The stringy cohomology groups
for toric DM stacks are particularly accessible. Like the cohomology and Grothendieck groups
of toric varieties, the cohomological invariants of toric DM stacks can be computed  and studied explicitly
in terms of the combinatorial data, the stacky fan, that is used to define the stack.  That combinatorial data is 
encoded in the homogeneous coordinate ring. 
 
  \subsection{}
 The author would like to thank Lev Borisov for pointing out that the Grothendieck group results in an earlier version of this paper applied only to complete toric stacks. That prompted the author to extend the result to
 other toric stacks. The author also thanks Paul Horja for pointing out a serious error in an earlier version of section \ref{sect.no.hypo}.

\section{The Grothendieck group of $\cX$}
\label{sect.K0}

 \subsection{}
 \label{sect.notn}
Let $G$ be a  closed subgroup of a torus  $(\GG_m)^r$ and write $\G := \Hom(G,\GG_m)$ 
for its rational character group.

Let $V$ be a finite dimensional rational representation of $G$. We introduce the following notation:
\begin{enumerate}
  \item 
 $\L$ is a set indexing a basis  $\{x_\rho \; | \; \rho \in \L\}$ for $V^*$ consisting of $G$-eigenvectors (in applications to a toric stack $\L$ will be the set of 1-dimensional rays in the fan);
  \item 
  $\L_1,\ldots,\L_n$ are subsets of $\L$; 
  \item 
  $L_m$, $1 \le m \le n$, is  the common zero locus of $\{x_\rho \; | \; \rho \in \L_m\}$;
  \item
   $Z:=L_1 \cup \ldots \cup L_n$;
   \item{}
   $\cX$ is the stack-theoretic quotient $[(V-Z)/G]$;
\item{}
For each $\a \in \G$,
\begin{enumerate}
  \item 
  $t^\a$  is the corresponding basis element of the integral
group ring $\ZZ\G$, and 
  \item 
   $\cO_{\cX}(\a)$ denotes the $\cO_\cX$-module that is $\cO_{V-Z}$ endowed with its
canonical $G$-equivariant structure twisted by the character  $\a$;
\end{enumerate} 
\item{}
 $S$:= the polynomial ring $k[x_\rho \; | \; \rho \in \L]=\cO(V)$ is endowed with the {\it right} $G$-action defined by $f^g(v):=f(gv)$ for $f \in S$, $g \in G$, and $v \in V$;
 \item{}
 $S$ is made into a $\G$-graded $k$-algebra with homogeneous components
$$
S_\a:=\{f \; | \; f^g =\a(g) f \; \hbox{for all $g \in G$}\}.
$$
\item{}
if $\rho \in \L$, we define the character $ |\rho|: G \to \GG_m$ by $(x_\rho)^g:=|\rho|(g) x_\rho$ for $g \in G$;
thus $x_\rho \in S_{|\rho|}$.
 \end{enumerate}

\begin{thm}
\label{thm.K0}
The map $\G \to \coh \cX$, $\a \mapsto [\cO_{\cX}(\a)]$, induces a ring isomorphism
$$
K_0(\coh \cX) \cong {{\ZZ\G}\over{(q_1,\ldots, q_n)}}
$$
where 
\begin{equation}
\label{defn.q}
q_m:=\prod_{\rho \in \L_m} (1-t^{|\rho|}).
\end{equation}
\end{thm}

\subsection{Temporary Hypothesis}

We first prove Theorem \ref{thm.K0} under the following equivalent hypotheses:
\begin{enumerate}
  \item 
  the only Zariski-closed  $G$-orbit in $V$ is $\{0\}$;
  \item 
  the invariants $\cO(V)^G$ consist of the constant functions;
  \item 
  $S_0=k$.
\end{enumerate} 
This hypothesis applies from now through section \ref{sect.the.proof}.
In Section \ref{sect.LB} is an example where the hypothesis fails but we show how to side-step the
problem for  that example.
Section \ref{sect.no.hypo} shows the hypothesis can be removed for all toric DM stacks.

\subsection{}

We write $\Gr(S,\G)$ for the category of $\G$-graded $\cO(V)$-modules with degree-preserving
homomorphisms, and $\sgr( S,\G)$ for the full subcategory of finitely generated modules. Vistoli's result  \cite[Example 7.21]{V} establishes  an equivalence 
$$
\Qcoh \cX \equiv  \hbox{$G$-equivariant quasi-coherent  $\cO_{V-Z}$-modules} \equiv {{\Gr( S,\G)}\over{\sT(\fa)  }}
$$
where $\sT(\fa)$ is the full subcategory consisting of direct limits of 
finitely generated modules whose support is contained in $Z$.

We write $\fa_m$ for the ideal in $S$ generated by $\{x_\rho \; | \; \rho \in \L_m\}$ and 
$$
\fa:= \bigcap_{m=1}^n \fa_m. 
$$
Each $\fa_m$ is a prime ideal, $\fa$ is a radical ideal, $L_m$ is equal to the zero locus of $\fa_m$,
and $Z$, which is defined in section \ref{sect.notn}(4), is equal to the zero locus of $\fa$. 
 
We define 
$$
q(t): = \prod_{\rho \in \Lambda} \big(1-t^{| \rho|} \big) \in \ZZ\G.
$$

\subsection{Hilbert series.}

The hypothesis that $\{0\}$ is the only $G$-stable closed subvariety of $V$ implies that $S_0=k$.
It follows that every homogeneous component of $S$ has finite dimension. To prove this, consider
the ideal $\fb$ generated by a homogeneous component $S_\b$. Then 
$\fb$ is generated by a finite number of elements $w_1,\ldots,w_r$ belonging to $S_\b$. 
If $w$ is any element in $S_\b$, then $w \in \fb$ so $w=w_1z_1+\cdots w_r z_r$ for some $z_i$s in $S$.
Taking the degree-zero component of each $z_i$, it follows that 
$w \in kw_1+\cdots +kw_r$. Hence $\dim S_\b <\infty$.

Let $M$ be a finitely generated $\G$-graded $S$-module. Then each homogeneous component $M_\a$,
$\a \in \G$,  has finite dimension so we may define the {\sf Hilbert series} of $M$ to be the formal expression 
$$
H_M(t):= \sum_{\a \in \G} (\dim_k M_\a) t^\a.  
$$
For the twist $M(\b)$ we have $H_{M(\b)} (t)= t^{-\b}.H_M(t)$.  
If $U$ and $V$ are $\G$-graded vector spaces so is their tensor product and 
$H_{U \otimes V}(t)=H_U(t).H_V(t)$ provided the homogeneous components of all these vector
spaces have finite dimension. 
 Since
$$
S \cong \bigotimes_{\rho \in \Lambda} \CC[x_\rho]
$$
it follows that
$$
H_S(t) = \prod_{\rho \in \Lambda} (1+t^{|\rho|} + t^{2|\rho|} + \cdots) =  q(t)^{-1}. 
$$

Since $M$ has a projective resolution in $\Gr(S,\G)$  involving only a finite number of direct sums
of twists $S(\a)$ of the free module $S$,
$$
H_M(t)=f(t) H_S(t)
$$ 
for some $f(t) \in \ZZ\G$.  Therefore $H_M(t)q \in \ZZ\G$ and $H_M(t)$ is not just a formal expression but 
a well-defined element of the localized group ring $\ZZ\G[q^{-1}]$.

\subsection{An aside on projectives}
\label{sect.projs}

The following lemma is ``well-known'' when the graded algebra is assumed to be connected.
However, the word ``connected''  is usually applied to algebras graded by a {\it free} abelian group and
are {\it not} making that restriction in this paper. We therefore include a proof. 

\begin{lem}
\label{lem.projs}
Let $k$ be a field and $(S,\G)$ a graded $k$-algebra such that $S_0=k$.
Suppose the only homogeneous units in $S$ are the elements in $k-\{0\}$. Then
\begin{enumerate}
  \item  
  $\fm:=\sum_{\a \ne 0} S_\a$ is the unique maximal   graded ideal of $S$;
  \item 
  if $M$ is a finitely generated graded $S$-module such that $\fm M=M$, then $M=0$;
  \item 
  if $P$ is a finitely generated projective graded $S$-module, then $P$ is isomorphic to a direct
  sum of graded free modules $S(\a)$ for various $\a$s in $\G$;
  \item{}
  two graded free modules $P:=S(\a_1) \oplus \cdots \oplus S(\a_m)$ and 
  $Q:=S(\b_1) \oplus \cdots \oplus S(\b_n)$ are isomorphic if and only if 
  there is an equality of multi-sets $\{\{\a_1,\ldots,\a_m\}\}=\{\{\b_1,\ldots,\b_n\}\}$. 
\end{enumerate}
\end{lem}
\begin{pf}
We do not assume $S$ is commutative in this proof. 

(1)
By hypothesis, $S=\fm \oplus k$ so, if $\fm$ failed to be a left ideal, $S\fm$ would contain $k$.
In particular, $S_\a \fm_{-\a}$ would be non-zero for some $\a \in \G$, so there would be elements $x \in S_\a$
and $y \in \fm_{-\a}$ such that $xy \ne 0$. Without loss of generality $xy=1$. Now $yx$ is also in $k$. It can't be zero because then $0=x(yx)=(xy)x=x$ which is absurd. Hence $x$ has a left and a right inverse, so 
$xy=yx=1$.  That contradicts the hypothesis about the homogeneous units in 
$S$ so we conclude that $S\fm=\fm$. By a similar argument $\fm$ is a right ideal. 

To see that $\fm$ is the unique maximal graded left ideal suppose $J$ is a graded left ideal
 that is not contained in $\fm$. Then $\fm +J=S$, from which it follows that $k \subset J$ and $J=S$. 

(2)
Suppose $\fm M=M$.
Because $M$ is finitely generated, if it were non-zero it would have a non-zero graded 
cyclic quotient module, $\Mol$ say. Suppose $\Mol$
is isomorphic to $S/J$ (with some shift in the grading). Then $J \subset \fm$ by (1), so $\fm(S/J) \ne S/J$. Hence $\fm \Mol\ne \Mol$. But this contradicts the fact that $\fm M=M$. We deduce that $M=0$.

(3)
Let $P$ be a non-zero finitely generated graded projective left $S$-module. Then $\fm P \ne P$.
Let $V$ be a graded subspace of $P$ such that $P=V \oplus \fm P$. Then $SV=P$ because 
$\fm(P/SV) =P/SV$. There is therefore a surjective 
degree-preserving homomorphism $\psi:S \otimes_k V \to P$, $\psi(s \otimes v)=sv$. 
Let $K=\ker \psi$. Because $P$ is projective applying $S/\fm \otimes_S -$
to the exact sequence $0 \to K \to S \otimes V \to P \to 0$ produces an exact 
sequence. It follows then that $K/\fm K=0$, so $K=0$ by (2)
and we deduce that $P \cong S \otimes V$, as required.

(4)
It follows from the argument in (3) that $P \cong Q$ if and only if the graded vector spaces $P/\fm P$ 
and $Q/\fm Q$ are isomorphic. But isomorphism of those two graded vector spaces is obviously equivalent to the condition
in (4).
\end{pf}

\subsubsection{Remarks}
\label{sect.remks.projs}
Given the result in Lemma \ref{lem.projs}, it might be sensible to say that a graded $k$-algebra $(A,\G)$ is 
{\sf connected} if $A_0=k$ and the only homogeneous units in $A$ are the elements in $k-\{0\}$.

Suppose $(A,\G)$ is noetherian and connected in this sense and has finite global dimension. 
Let $\sT$ be any localizing subcategory of $\Gr(A,\G)$. Let $\cF$ be the image of a finitely generated graded $A$-module in $\Gr(A,\G)/\sT$. Then $\cF$ has a finite resolution in  $\Gr(A,\G)/\sT$ in which each term is a direct sum of various twists of $\cO$, where $\cO$ denotes the image of $A$ in  $\Gr(A,\G)/\sT$. 

This remark is applied to certain DM stacks in Proposition \ref{prop.res} in  section \ref{ssect.hcrs}.

\subsection{}
We now return to the main line of the proof, so $S$ once more denotes
 the polynomial ring $\cO(V)$ which is assumed to satisfy the 
hypotheses stated at the beginning of section \ref{sect.K0}. 

\begin{prop}
\label{prop.K0}
There are mutually inverse $\ZZ$-algebra isomorphisms
\begin{equation}
\label{map.K0}
K_0(\sgr S) \to \ZZ\G, \qquad [M] \mapsto H_M(t) q
\end{equation}
and 
$$
\ZZ\G \to K_0(\sgr S), \qquad t^\a \mapsto [S(-\a)].
$$
\end{prop}
\begin{pf}
Let $\sP$ denote the full subcategory of $\sgr(S,\G)$ consisting of the projective modules.
By Lemma \ref{lem.projs},  every finitely generated projective graded  $S$-module is isomorphic to 
a unique finite direct sum of various 
$S(\a)$s, so the map $t^\a \mapsto [S(-\a)]$  is an isomorphism from $\ZZ\G$ to $K_0(\sP)$.
This is an isomorphism of rings because $S(\a) \otimes_S S(\b) \cong S(\a+\b)$. 

Every $M \in \sgr (S,\G)$ has a finite resolution by finitely generated projective graded $S$-modules so
the inclusion functor $\sP \to \sgr (S,\G)$ induces an isomorphism of Grothendieck groups 
$K_0(\sP) \to K_0(\sgr S).$ We therefore have a ring isomorphism 
$$
\Psi:\ZZ\Gamma \to K_0(\sgr S), \qquad \Psi(t^\a):= [S(-\a)].
$$

If $0 \to L \to M \to N \to 0$ is exact in $\sgr S$, then $H_M=H_N+H_L$, so
the universal property of $K_0$ ensures here is a well-defined group homomorphism
$$
\Phi:K_0(\sgr S) \to \ZZ \G, \qquad \Phi([M])=H_M(t)q.
$$
In particular, $\Phi([S(\a)]) = H_{S(\a)}(t)q=t^{-\a} = \Psi^{-1}([S(\a)])$. Since $\ZZ\G$ is spanned by the 
$t^\a$s, $\Phi=\Psi^{-1}$. 
\end{pf}

For an arbitrary pair of modules $M, N \in \sgr S$, the usual
argument shows that
$
[M].[N] = \sum_{i \ge 0} (-1)^i \Tor^S_i(M,N)
$ 
where the Tor-groups are computed as graded $S$-modules.

\subsection{Proof of Theorem \ref{thm.K0}}
\label{sect.the.proof}
Let $\sT=\sT(\fa) \cap \sgr S$. Thus $\sT$ is the full
 subcategory of $\sgr S$ consisting of those modules supported on $Z$, i.e., the finitely generated graded 
  modules annihilated by a suitably large power  of $\fa$. 
 
The localization sequence for $K$-theory gives an exact sequence
\begin{equation}
\label{loc.seq}
K_0(\sT)  \stackrel{\iota}{\longrightarrow}   K_0(\sgr S) \longrightarrow K_0(\coh \cX)   \longrightarrow 0.
\end{equation}
If $M \in \sT$ so is $M(\a)$ for all $\a \in \G$ so $K_0(\sT)$ is a $\ZZ\G$-module under the action
$t^\a.[M]=[M(-\a)]$. The arrow $\iota$ in (\ref{loc.seq}) is induced by the inclusion 
$\sT \to \gr  S$ so is a $\ZZ\G$-module homomorphism. Therefore, after identifying $K_0(\sgr S)$ with $\ZZ\G$ as in Proposition \ref{prop.K0}, $K_0(\coh\cX)$ is isomorphic to $\ZZ\G$ modulo the ideal generated by the images of a set of $\ZZ\G$-module generators for $K_0(\sT)$. 

By definition, $\sT$ consists of the modules annihilated by a power of $\fa$ so, by d\'evissage, the natural
map $K_0(\sgr S/\fa) \to K_0(\sT)$ is an isomorphism, even an isomorphism of $\ZZ\G$-modules. 
Since ${\fa}$ is the intersection of the $\fa_m$s every $M \in \gr(S/\fa)$ has a finite filtration $M=M_0 \supset M_1 \supset \dots \supset M_r=0$ such that each
slice $M_i/M_{i+1}$ a finitely generated graded $S/\fa_m$-module for some $m$. 
Since $S/\fa_m$ is a polynomial
ring, $M_i/M_{i+1}$ has a finite resolution as an $S/\fa_m$-module in which
all the terms are direct sums of various twists $(S/\fa_m)(\a)$. It follows that $K_0(\sgr S/\fa)$, and hence 
$K_0(\sT)$, is generated as a $\ZZ\G$-module by the classes $[S/\fa_m]$, $1 \le m  \le n$.

The image of $[S/\fa_m]$ under the first map in (\ref{loc.seq}) is $[S/\fa_m]$. Since $S/\fa_m$ is the polynomial ring on the indeterminates $\{x_\rho \; | \; \rho \in \L-\L_m\}$, 
$$
H_{S/\fa_m}(t)=  \prod_{\rho \in\L- \L_m} (1-t^{|\rho|})^{-1}.
$$
Under the isomorphism in (\ref{map.K0}), the image of $[S/\fa_m]$ in $\ZZ\G$ is therefore
$$
H_{S/\fa_m}(t)q= q \prod_{\rho \in\L- \L_m} (1-t^{|\rho|})^{-1} =   \prod_{\rho \in  \L_m} (1-t^{|\rho|}).
$$
This completes the proof of Theorem \ref{thm.K0}.  \hfill $\square$

\subsection{Example}
\label{sect.LB}

This example was prompted by a question of Lev Borisov.  

Let $Bl_{(0,0)}\CC^2$ denote the blowup of $\CC^2$ at the origin. The usual fan for this toric variety is 
that spanned by  $(1,0)$, $(1,1)$, and $(0,1)$. Cox's homogeneous coordinate ring
is $(\CC[x_0,x_1,x_2],\ZZ,(x_0,x_2))$ where the grading is given by $\deg x_i=(-1)^i$. Since $\CC[x_0,x_1,x_2]_0 \ne \CC$, Theorem \ref{thm.K0} does not apply. 

However, Theorem  \ref{thm.K0} {\it does} apply if we present $Bl_{(0,0)}\CC^2$ as an open subscheme of  the Hirzebruch surface $\FF_1$ because then $Bl_{(0,0)}\CC^2$ has homogeneous
coordinate ring $\CC[t_0,t_1,x_0,x_1]$ with $\ZZ^2$-grading $\deg t_0=\deg t_1=(1,0)$, $\deg x_0=(-1,1)$,
and $\deg x_1 = (0,1)$ with irrelevant ideal $x_1(t_0,t_1)$, and $\CC[t_0,t_1,x_0,x_1]_{(0,0)}=\CC$.
The locus $Z$ is the union of the subspaces
$x_1=0$ and $t_0=t_1=0$. The group algebra $\ZZ\G$ is $\ZZ[u^{\pm 1},v^{\pm 1}]$
with $(1,0)=u$ and $(0,1)=v$ so
$$
K_0( Bl_{(0,0)}\CC^2) = {{\ZZ[u^{\pm 1},v^{\pm 1}]}\over{(1-v,(1-u)^2)}}.
$$
For the first presentation of $Bl_{(0,0)}\CC^2$  as $\CC^3-Z(x_0,x_2)/\CC^\times$ there are non-trivial
closed orbits such as $x_0x_1=x_2x_1=1$. However, for the second presentation 
 of $Bl_{(0,0)}\CC^2$  as $(\CC^4-Z(x_1t_0,x_1t_1))/\CC^\times \times \CC^\times$ the
 only closed orbit is the origin.

\subsection{Removing the temporary hypothesis}
\label{sect.no.hypo}

Let $\cX=[(V-Z)/G]$ be as described at the beginning of section \ref{sect.notn}, but {\it do not assume}
that $\{0\}$ is the only closed $G$-orbit. We will show there is alternative data $V'$, $Z'$, $G'$ such
that 
\begin{enumerate}
  \item 
  the stack $\cX':=[(V'-Z')/G']$ is isomorphic to $\cX$ and
  \item  $\{0\}$ is the only closed orbit  for the $G'$ action on $V'$. 
\end{enumerate} 
We will do this by showing that the data $(S,\G,\fa)$, which determines and is determined by $(V,G,Z)$,
 may be replaced by data $(S',\G',\fa')$ such that the degree zero component of $S'$ is $k$ and 
$$
 {{\Gr(S,\G)}\over{\sT(\fa)}} \equiv  {{\Gr(S',\G')}\over{\sT(\fa')}}.
 $$

Let $(A,\Ups)$ be a graded ring. If $\Upsol$ is a quotient of $\Ups$, then $A$ becomes a $\Upsol$-graded
ring with respect to the grading
$$
A_\d := \sum_{i \in \d} A_i
$$
for each coset $\d \in \Upsol$. An ideal $I$ in $A$ that is not graded for the $\Ups$-grading might be
graded with respect to the $\Upsol$-grading, and in that case the natural  homomorphisms to the quotients 
give a morphism  of  graded rings $(\psi,\theta):(A,\Ups) \to (A/I,\Upsol)$. This idea is used in the next result.

\begin{prop}
\label{xxprop1.9.1}
\cite[Prop. 1.9.1]{SZ},
Let $(A,\Ups)$ be a graded ring having central homogeneous units $z_1,\ldots,z_d$ and define
$$
\Aol:={{A}\over{(z_1-1,\ldots,z_d-1)}}.
$$
Suppose that the subgroup  $\Sigma$  generated by $\{\deg z_i \; | \; 1 \le i \le d\}$ is free of rank $d$.
If there is a subgroup $\G\subset \Ups$ such that $\Ups=\G\oplus \Sigma$, 
there is an equivalence of categories
$$
\UseComputerModernTips
  \xymatrix{
 F^*:   \Gr (A,\Ups) \ar[rr]^>>>>>>>>{\sim}  && \Gr (\Aol,\Ups/\Sigma) , \qquad F^*M:=\Aol \otimes_A M.
}
$$
 \end{prop}

\begin{cor}
\label{cor.BwtB}
Let $(S,\G)=k[x_\rho \; | \; \rho \in \L]$ be a graded polynomial ring as in section \ref{sect.notn}.
Let $$\wtS=S\otimes_k k[z^{\pm 1}]$$ where the $z$ is a central indeterminate  
and give $\wtS$ a $\G \times \ZZ$ grading by declaring that
$$
\deg z :=(0,1)
\qquad \hbox{and} \qquad \deg_{\wtS} x_\rho :=(\deg_Sx_\rho,1).
$$
Then
 $\Gr(S,\G) \equiv \Gr(\wtS,\G \times \ZZ)$.  
\end{cor}
\begin{pf}
This follows from Proposition \ref{xxprop1.9.1} with 
$\Ups = \G \times \ZZ$, $\Sigma=(0,\ZZ)$,  $\D=(\G,0)$, $\wtS$ playing the role of $A$. and $S$ playing the role of $\Aol$.
\end{pf}

\begin{prop}
\label{prop.connected}
Retain the notation in Corollary \ref{cor.BwtB}. Let $S'$ be the subring $S[z]$ of $\wtS$. 
Suppose that $\fa$ is a graded ideal in $(S,\G)$ and let $\fa':=S'\fa z$.
Then 
$$
{{\Gr(S,\G)}\over{\sT(\fa)}} \equiv {{\Gr(S',\G\times \ZZ)}\over{\sT(\fa')}} 
$$
Furthermore, $S'_{(0,0)}=k$. 
\end{prop}
\begin{pf}
For any graded ring $(S',\Ups,\fa')$ having a central homogeneous regular element $z$ such that 
$\fa' \subset zS'$ the induction functor $S'[z^{-1}] \otimes_{S'} -$ from $\Gr(S',\Ups)$ to
$ \Gr(S'[z^{-1}],\Ups)$ induces an equivalence  
$$
 {{\Gr(S',\Ups)}\over{\sT(\fa')}} \equiv {{\Gr(S'[z^{-1}],\Ups)}\over{\sT(\fa'[z^{-1}])}}.
 $$
 (This is an analogue of the fact that if $R$ is a commutative ring with an ideal $\fa$
 contained in a principal ideal $zR$, then $\Spec R-Z(\fa) = \Spec R[z^{-1}] - Z(\fa[z^{-1}])$.)
 
 Applying this to the case of interest with $\Ups=\G \times \ZZ$,
$$
 {{\Gr(S',\Ups)}\over{\sT(\fa')}} \equiv {{\Gr(S'[z^{-1}],\Ups)}\over{\sT(\fa'[z^{-1}])}}
 = {{\Gr(\wtS,\Ups)}\over{\sT(\wtS \fa')}} = {{\Gr(\wtS,\Ups)}\over{\sT(\wtS \fa)}} .
 $$
Under the equivalence in Corollary \ref{cor.BwtB}, the graded $S$-modules annihilated by 
 a power of $\fa$ correspond to the graded $\wtS$-modules annihilated by 
 a power of $\wtS\fa$, so the equivalence in Corollary \ref{cor.BwtB} induces an equivalence
 $$
{{\Gr(S,\G)}\over{\sT(\fa)}} \equiv   {{\Gr(\wtS,\Ups)}\over{\sT(\wtS \fa)}} 
$$
 between the quotient categories. This completes the proof of the claimed equivalence of categories.

The degree zero component of $S'$ with its $\G \times \ZZ$-grading
is spanned by the homogeneous 
elements $x z^{s}$ such that $x$ is a word of length $r$ in the letters $x_\rho$ and  
$$
\deg_{\wtS}(xz^s) = (\deg_S x, r)+(0,s)  =(0,0) \in \G \times \ZZ.
$$
It follows that $S'_{(0,0)}=k$.
\end{pf}

We now define $V'=\Spec S'$, $Z'=Z(\fa')$, and $G'=\Spec k \G'$ where $\G'=\G \times \ZZ$.  
Because $S'=S \otimes_k k[z] $, $V'=V \times k$. Because $\fa'=S'\fa z$,
$$
Z'= (Z \times k) \cup (V \times \{0\})  
$$
We write the group ring for $\G'$ as 
$$
\ZZ\G'=\ZZ\G[t^{\pm 1}].
$$
where $t=\deg_{S'}z$.
Because the degree zero component of $S'$ is $k$ Theorem \ref{thm.K0} gives
$$
K_0(\cX') \cong {{\ZZ\G'}\over{(q_1,\ldots,q_m,q)}}
$$
where $q_1,\ldots,q_m$ have the same meaning as before, and $q=1-t$. 
Therefore
$$
K_0(\cX') \cong {{\ZZ\G}\over{(q_1,\ldots,q_m)}}.
$$

The equivalence of categories in Proposition \ref{prop.connected} says that $\cX\cong \cX'$, but one can also 
see this geometrically because 
$$
\biggl[{{V'-Z'}\over{G'}}\biggr]=\biggl[{{(V-Z) \times (k-\{0\})}\over{G\times \GG_m}}\biggr].
$$
Let $\eta \in \GG_m$. 
The $\G\times \ZZ$-grading on $S'$ defined in Corollary \ref{cor.BwtB} is such that  the action of $(1,\eta)\in G'$ on a point in $(v,\l) \in V'= V \times k$  is given by 
$$
(1,\eta).(v,\l) = (\eta v, \eta \l).
$$
It is now clear that the origin of $V'$ is in the closure of every $G'$-orbit on $V'$.

 \section{Homogeneous coordinate rings for some stacks}\footnote{The idea of using $(\cO(V), \G)$ as a ``homogeneous coordinate ring'' of a stack is developed more fully in
\cite{SZ} although the main focus there is on homogeneous coordinate rings of non-commutative schemes.}
\label{ssect.hcrs}

Let $k$ be a field, $\cX$ a stack over $\Spec k$, and suppose we have data $(A,\D,\fa)$ consisting of 
\begin{enumerate}
  \item 
  an abelian group $\D$, 
  \item 
  a $\D$-graded commutative $k$-algebra $A$, and
  \item{}
   a graded ideal $\fa$. 
\end{enumerate}
Let $G$ be the affine group scheme $\Spec k\D$ where $k\D$ is given its natural Hopf algebra structure.
Let $Z(\fa)$ denote the zero locus of $\fa$. 
We call $(A,\D,\fa)$, or simply $A$ if the other data is clear from the context,  a {\sf homogeneous coordinate ring} of $\cX$ if 
$$
\cX \cong  \biggl[{{\Spec A -Z(\fa)}\over{G}}\biggr].
$$
If $(A,\D,\fa)$ is a homogeneous coordinate ring for $\cX$, Vistoli's result  \cite[Example 7.21]{V}  
tells us there is an equivalence of monoidal categories 
\begin{equation}
\label{eq.qchX}
\Qcoh \cX \equiv {{\Gr(A,\D)}\over{\sT(\fa)}}
\end{equation}
where $\sT(\fa)$ is  the localizing subcategory consisting of the graded modules $M$ such that
$H_\fa^0(M)=M$.

Let 
$$
\pi^*:\Gr(A,\D) \to \Qcoh \cX
$$
be the functor inducing the equivalence in (\ref{eq.qchX}).  The functor $\pi^*$ is analogous to the functor
$$
M \rightsquigarrow \wtM
$$
that is used in the classical case for schemes of the form $\Proj A$ where $A$ is an $\NN$-graded 
commutative ring generated as an $A_0$-algebra by $A_1$. Indeed, when $A$ satisfies those hypotheses
and $\fa=A_{\ge 1}$, then $\pi^*$ is  the functor $M \rightsquigarrow \wtM$.

Because $\cO_\cX$ and $\pi^* A$ are neutral objects for the internal tensor product on the two categories 
in (\ref{eq.qchX}) we can replace $\pi^*$ by its composition with a 
suitable auto-equivalence of  $\Qcoh \cX$  and so assume that $\pi^* A = \cO_\cX$. 
We will assume this has been done. 

The general results on quotient categories  in \cite{Gab} tell us that $\pi^*$ is exact and has a right  adjoint
that we will denote by $\pi_*$.  Furthermore,  the counit is an isomorphism 
$\pi^*\pi_* \cong \id_{\Qcoh \cX}$ and the unit fits into an exact sequence
$$
0 \to  H^0_\fa(M) \to M \to \pi_*\pi^*M \to H^1_\fa(M) \to 0
$$
that is functorial in $M$.

The following conditions on a graded $A$-module  $M$ are equivalent: 
\begin{enumerate}
\item
$M \in \sT(\fa)$;
  \item 
  $\pi^*M=0$;
  \item 
  $H^0_\fa(M)=M$;
  \item 
  the support of  every finitely generated submodule of $M$ is contained in $Z(\fa)$.
\end{enumerate}
We call $M$ a {\sf torsion} module if it satisfies (1)--(4). Always $H_\fa^0(M)$ is the largest submodule of $M$ that $\pi^*$ sends to zero. We say $M$ is {\sf torsion-free} if 
$H_\fa^0(M)=0$. 

\subsubsection{Resolutions by direct sums of invertible $\cO_\cX$-modules}

The following is an immediate consequence of Lemma \ref{lem.projs} and the remarks in 
section \ref{sect.remks.projs}.

\begin{prop}
\label{prop.res}
Let $\cX$ be a stack and  $\cF \in \coh \cX$. Suppose that $\cX$ 
has a homogeneous coordinate ring $(A,\D,\fa)$  such that $A$ noetherian, has finite global dimension, and
 is connected in the sense of the remark in section \ref{sect.projs}. Then $\cF$ has a finite resolution in 
 which every term is a direct sum of invertible $\cO_\cX$-modules of the form $\cO_\cX(\a)$ for various 
 $\a$s in $\D$.
\end{prop}

\section{The  Picard group when $\cX$ has a homogeneous coordinate ring}
\label{sect.picard}

In this section we compute the Picard group $\Pic \cX$ which is, 
by definition,  the group of isomorphism classes of  invertible $\cO_\cX$-modules with group 
operation given by $\otimes$.

\subsection{Graded domains}

 Suppose $(A,\D)$ is a {\sf graded-domain}, i.e.,  every non-zero homogeneous element 
is regular, i.e., not a  zero-divisor. Then $(A,\D)$ embeds in its graded ring of fractions
$$
K:=\{ab^{-1} \; | \; \hbox{$a$ and $b$ are homogeneous and $b\ne 0$}\}
$$
with grading given by $\deg(ab^{-1})=\deg a - \deg b$.  

Every $\D$-graded $K$-module is isomorphic to a direct sum of various twists $K(\a)$, $\a \in \D$. Furthermore, $K(\a) \cong K$  if and only if $K_\a \ne 0$.\footnote{We will not need this fact, but 
 $\Gr(K,\D)$ is equivalent to the category of locally unital modules over  the 
direct sum $K_0^{\oplus |\D/\G|}$ where $\G=\{\a \in \D \; | \; K_\a \ne 0\}.$
The subgroup $\G$ is  generated by $\{\a \; | \; A_\a \ne 0\}$.} 
  
 A non-zero homogeneous non-unit $a \in A$ is said to be {\sf graded-irreducible} if in every factorization $a=bc$ in which $b$ and $c$ are homogeneous either $b$ or $c$ is a unit. A non-zero homogeneous non-unit $a \in A$ is said to be {\sf graded-prime} if 
whenever $a$ divides a product $bc$ of homogeneous elements $a$ divides either $b$ or $c$.
We say that $(A,\D)$ is {\sf graded-factorial} if 
every homogeneous element is a product of graded-prime elements.
If $A$ is graded factorial and noetherian, 
then every non-zero homogeneous element is either a unit or a product of
graded-irreducible elements in a unique way; the notion of {\it greatest common homogeneous 
divisor} for a set of homogeneous elements of $A$ therefore makes sense.

We say $A$ is {\sf graded-noetherian} if every graded ideal of $A$ is finitely generated.

\subsection{Invertible $\cO_{\cX}$-modules} 
\label{sect.pic.pfs}

In all the results in this section we assume  that $(A,\D,\fa)$ satisfies conditions (1)-(3) at the 
beginning of section \ref{ssect.hcrs} and that it is a homogeneous coordinate ring  for a 
stack $\cX$.

 \begin{lem}
 \label{lem.pic}
 Suppose $(A,\D,\fa)$ is a homogeneous coordinate ring for $\cX$.
 Suppose further that $A$ is graded-noetherian graded-factorial graded-domain such that 
 $H^0_\fa(A) = H^1_\fa(A) =0$.
Let $M$ and $N$ be finitely generated graded $A$-modules, and suppose there is a degree-preserving 
homomorphism
$\phi:M \otimes_{A} N \to A$ such that $\pi^*\phi$ is an isomorphism. 
Then $\pi^*M \cong \cO_\cX(\a)$ and $\pi^*N \cong \cO_\cX(-\a)$ for some $\a \in \D$.
 \end{lem}
\begin{rem}
The hypothesis $H^0_\fa(A) = H^1_\fa(A) =0$ implies that the natural map 
$A(\a) \to \pi_*\pi^* A(\a) =  \pi_* \bigl(\cO_\cX(\a)\bigr)$ is an isomorphism for all $\a \in \D$. We will  write $A(\a) 
=  \pi_*\cO_\cX(\a)$ to denote this fact.
\end{rem}
 \begin{pf}
Because $\pi^*\phi$ is an isomorphism, $\ker \phi$ and $\coker\phi$ are both torsion.
  
First we prove the result under the assumption that $H_\fa^0(M) = H_\fa^0(N)=0$.

The result is vacuous if $\fa=0$, so we assume $\fa \ne 0$. Every non-zero homogeneous 
element in $K$ is a unit so  $\fa K=K$. It follows that the functor
$K \otimes_A -$ sends all torsion modules to zero. In particular, $K \otimes_A -$ kills the kernel and cokernel of $\phi$. But $K \otimes_A -$ is an exact functor so $\phi$  induces an isomorphism $(M \otimes_A K) \otimes_K (N \otimes_A K) \to K$. Since every homogeneous element in $K-\{0\}$ is a unit,  
 there is an element $\a \in \D$ such that $M \otimes_A K \cong K(\a)$ and $N \otimes_A K \cong K(-\a)$. Let $f$ and $g$ be the obvious compositions $M \to M \otimes_A K \to K(\a)$ and $N \to N \otimes_A K \to K(-\a)$ respectively. 
 
 Let $\mu$ be the restriction to $f M \otimes_A g N$ of the multiplication map $K(\a) \otimes_A K(-\a) \to K$ and consider the not-necessarily-commutative diagram
 $$
 \begin{CD}
 M \otimes_A N @>{\phi}>> A
 \\
 @V{f \otimes g}VV @VV{\iota}V
 \\
 f M \otimes_A g N @>>{\mu}> K
 \end{CD}
 $$
 where $\iota:A \to K$ is the inclusion. 
 Since $K \otimes_A M \otimes_A N \cong K$, $\Hom_A(M \otimes_A N,K) \cong K$.  
 Hence there is some $c \in K_0$ such that $c\iota \circ \phi=\mu \circ (f
 \otimes g)$. Now replace $f$ by $c^{-1} f$ so that $\iota \circ \phi=\mu \circ (f \otimes g)$, and replace $M$ and $N$ by $f M$ and $g N$, so that $M$ and $N$ are graded $A$-submodules of $K(\a)$ and $K(-\a)$ 
 respectively, and $\phi$ is the restriction of the multiplication map.
 The image of $\phi$ is therefore the product $MN$.   
 but $\pi^*(\coker\phi)=0$ so $A/MN$ is torsion.  
 
Because $MN \subset A$, there is a non-zero homogeneous element $a \in A$ such that 
$Ma \subset A$. 
 Hence 
$Ma=\sum_{j=1}^n Ab_j$ for some  homogeneous elements $b_j \in A$, $1 \le j \le n$. 
 Let $d$ be the greatest common homogeneous divisor of the $b_j$s. Then $\{q \in K \; | \; Maq \subset A\} = Ad^{-1}$. 
 Therefore  $a^{-1}N  \subset Ad^{-1}$ and  
  $$
  MN = Maa^{-1}N \subset Mad^{-1}A \subset A.
  $$
   It follows that 
 $A/ad^{-1}M$ is torsion, whence  $\pi^*M \cong \cO(\deg a - \deg d)$.
 
 This completes the proof when $M$ and $N$ are torsion-free. Now we deal with the general case.
 The map $\phi:M \otimes_A N \to A$ factors as a composition
  $$
  \begin{CD}
 M \otimes_A N @>{\phi_1}>> {{M}\over{\tau M}}  \otimes_A {{N}\over{\tau N}}   @>{\phi_2}>> A
 \end{CD}
 $$
 Since $\pi^*\phi$  is an isomorphism, $\pi^*\phi_1$ is monic; but $\phi_1$ is epic so $\pi^*\phi_1$ is epic too;
 hence $\pi^*\phi_1$, and therefore $\pi^*\phi_2$, is an isomorphism. 
 By the first part of the proof applied to $\phi_2$,   $\pi^*(M/\tau M) \cong \cO(\a)$ for some 
 $\a \in \D$. But $\pi^*M \cong \pi^*(M/\tau M)$, so  $\pi^*M \cong \cO(\a)$.
  \end{pf}

\begin{lem}
\label{lem.twists}
Suppose $(A,\D,\fa)$ is a homogeneous coordinate ring for the stack $\cX$.
 Suppose further that $A$ is graded-noetherian graded-factorial graded-domain such that 
 $H^0_\fa(A) = H^1_\fa(A) =0$. Then 
$\cO_\cX(\a) \cong \cO_\cX(\b)$ if and only if $A_{\a-\b}$ contains a unit.
\end{lem}
\begin{pf}
Multiplication by a unit $u \in A_{\a-\b}$ produces an isomorphism $g:A(\b) \to A(\a)$
of graded $A$-modules. Applying $\pi^*$ to $g$ produces an isomorphism 
$\pi^*g:\cO(\b)  \to \cO(\a)$.

Conversely, suppose that $f:\cO(\b) \to \cO(\a)$ is an isomorphism. Since $\pi^*\pi_* \cong \id$, 
applying $\pi^*$ to the map  $\pi_*f: \pi_*\cO(\b) \to \pi_*\cO(\a)$ produces $f$ again. 
The kernel and cokernel of $\pi_*f$ are therefore
torsion. But  
 $\pi_*\cO(\b) = A(\b)$, so $\ker( \pi_*f) =0$. However, $H^0_{\fa}(\coker(\pi_*f )  ) = \coker (\pi_*f)$ and 
 $H^1_\fa(A(\b)=0$ so  $\Ext^1_A(\coker(\pi_*f),A(\b))=0$. Consequently the  exact sequence
$$
\begin{CD}
0 @>>> A(\b) @>{\pi_*f}>> A(\a) @>>> \coker(\pi_*f) @>>> 0
\end{CD}
$$
splits. But $H^0_\fa(A(\a))=0$ so we conclude that $\coker(\pi_*f)=0$. Hence $\pi_*f$ is an isomorphism.
But every $\Delta$-graded $A$-module homomorphism $A(\b) \to A(\a)$ is multiplication 
by an element of $A_{\a-\b}$, so the required unit  exists.
\end{pf}

\begin{thm}
\label{thm.pic}
Suppose $(A,\D,\fa)$ is a homogeneous coordinate ring for the stack $\cX$.
 Suppose further that $A$ is graded-noetherian graded-factorial graded-domain such that 
 $H^0_\fa(A) = H^1_\fa(A) =0$.  
Define $$
\D_u:=\langle \a \in \D \; | \; A_\a 
\hbox{ contains a unit}\rangle.
$$
The map $\a \mapsto \cO(\a)$ induces an isomorphism 
$$
\begin{CD}
\D/ \D_u @>{\sim}>> \Pic(\cX).
\end{CD}
$$
\end{thm}
\begin{pf}
If $M \in \Gr(A,\D)$ is invertible so is $\pi^*M$, so the rule $\a \mapsto \cO(\a)$ is a  
homomorphism $\D \to \Pic \cX$. By Lemma \ref{lem.twists}, the kernel of this map
 is $\D_u$. 
It remains to show that the $\cO(\a)$s are the only invertible $\cO_{\cX}$-modules up to isomorphism.

 To this end, suppose $\cM \otimes \cN\cong\cO_\cX$. 
 Suppose $\cM=\pi^* M$ and $\cN=\pi^* N$.  By adjointness,
 the isomorphism $\cM \otimes \cN \to \cO_\cX$ is induced by a map $\phi:M \otimes_A N \to \pi_*\pi^* A$ whose  kernel and cokernel are torsion. But   $\pi_* \pi^*  A= A$ 
 so $\phi:M \otimes_A N \to A$.  The result now follows from   Lemma \ref{lem.pic}.
   \end{pf}

\begin{prop}
\label{prop.Pic.open}
Suppose $(A,\D,\fa)$ is a homogeneous coordinate ring for the stack $\cX$.
 Suppose further that $A$ is graded-noetherian graded-factorial graded-domain such that 
 $H^0_\fa(A) = H^1_\fa(A) =0$.
Let $f$ be a homogeneous element in $A$ of degree $\a$ and let $\cZ \subset \cX$
be the zero locus of $f$. 
Then there is an exact sequence
$$
\begin{CD}
\ZZ @>>> \Pic \cX @>>> \Pic(\cX -\cZ) @>>> 0
\end{CD}
$$
where the first map is $n \mapsto \cO(n\a)$.
\end{prop}
\begin{pf}
Write $\cU=\cX-\cZ$ and let $\iota:\cU \to \cX$ be the inclusion. Then 
$$
\Qcoh\, \cU \equiv{{\Gr(A[f^{-1}],\D)}\over{\sT'}}
$$
where $\sT'$  is the full subcategory consisting of $A[f^{-1}]$-modules with the property that
the support of all their finitely generated submodules is contained in $Z(\fa[f^{-1}])$. 

The inclusion $\iota$ is induced by the natural map $A \to A[f^{-1}]$.

But $A[f^{-1}]$  is graded factorial because $A$ is,  so there is a commutative diagram $$
\begin{CD}
\D/\D_u @>>> \D/\langle \D_u, \a  \rangle
\\
@VVV @VVV
\\
\Pic \cX @>>{\iota^*}> \Pic \cU
\end{CD}
$$ 
 in which the vertical maps are isomorphisms.
Since the upper horizontal arrow is surjective so is the lower one. The result now follows from
the fact that kernel of the upper arrow is  the map $\ZZ \to \D/\D_u$, $1\mapsto \a$.
\end{pf}

\section{Some examples}

\subsection{Stacky weighted projective spaces}
\label{sect.stacky.wPn}

Let $Q=(q_0,\ldots,q_n)$ be a sequence of positive integers. The weighted projective space $\PP(Q)$ is 
the scheme $\Proj A$ where $A$ is the weighted polynomial ring
$$
A=k[x_0,\ldots,x_n], \qquad \deg x_i=q_i.
$$
If the characteristic of $k$ does not divide any of the $q_i$s, then $\PP(Q)$ can be expressed 
as the quotient of $\PP^n$ modulo the coordinate-wise action of $$\mu_Q:=\mu_{q_0} \times \cdots \times 
\mu_{q_n},$$  where $\mu_q$ denotes the group of $q^{\th}$ roots of 1 in $k^\times$.

The {\sf stack-theoretic weighted projective space}
$$
\PP[Q] = \PP[q_0,\ldots,q_n]
$$
is defined to be the stack-theoretic quotient 
$$
\biggl[{{\AA^{n+1} - \{0\}}\over{\GG_m}}\biggr]
$$
where $\xi \in \GG_m$ acts by 
$$
\xi\cdot (x_0,\ldots,x_n) = (\xi^{q_0}x_0,\ldots, \xi^{q_n} x_n).
$$

Let $\fm=(x_0,\ldots,x_n)$.   Then $(A,\ZZ,\fm)$ is 
a homogeneous coordinate ring for $\PP[Q]$. 

 \begin{thm}
 \label{thm.pic2}
 Suppose $n \ge 1$. With the above notation, 
 $$
 \Pic  \PP[q_0,\ldots,q_n] \cong \ZZ.
 $$
  \end{thm}
 \begin{pf}
 Because $A$ is a UFD in the usual sense it is graded factorial. Because $n \ge 1$, $H^1_\fm(A)=0$
 so the result follows from Theorem \ref{thm.pic}
 \end{pf}

 \subsubsection{}
 When $n=0$,   $\PP[q]$ is the classifying stack $B(\mu_q)$. 
 In that case, rather than using $(k[x],\ZZ,(x))$ as a homogeneous coordinate ring for $\PP[q]$, we may 
 use  $(k[x,x^{-1}],\ZZ,0)$ as the homogeneous coordinate ring. 
 The  hypotheses of Theorem \ref{thm.pic}  are satisfied by  $(k[x,x^{-1}],\ZZ,0)$ 
 and $k[x,x^{-1}]$ has a unit of degree $nq$ for all $n$,  so $\Pic \PP[q]\cong \ZZ/q\ZZ$.
 
 \subsubsection{}

Let $\cM_{1,1}$ be the fine moduli space of pointed elliptic curves over $\CC$, and $\cMol_{1,1}$ its
usual compactification. It is well-known that $\cMol_{1,1} \cong \PP[4,6]$. Because $\cM_{1,1}=
\cMol_{1,1}-\{p\}$ where $p$ is the zero locus of a degree 12 element, Proposition \ref{prop.Pic.open}
gives 
$$
\Pic \cM_{1,1} \cong {{\ZZ}\over{12\ZZ}}.
$$

\subsection{Rugby balls}

Fix positive integers $p$ and $q$.
The orbifold whose underlying manifold is the Riemann sphere endowed with
the  groupoid structure given by cyclic groups of orders $p$ and $q$ at the north and south poles is called a {\sf rugby ball}.  If $p$ or $q$ is 1, it is called a {\sf teardrop}, sometimes {\sf Thurston's teardrop}.

Let 
$$
G:=  \bigl\{(\l_1,\l_2) \in \CC^\times \times \CC^\times \; | \; \l_1^p=\l_2^q\bigr\}
$$
 act on $\CC^2$ by component-wise multiplication and define  the  {\sf $(p,q)$-rugby ball} to be the stack
 $$
\FF[p,q]:=  \biggl[{{\CC^2-\{0\}}\over{G}} \biggr].
$$
 Give  the polynomial ring $S:=\CC[x,y]$ a grading by  the group 
 $\Gamma:=\langle e,e' \; | \; pe=qe' \rangle$ by declaring $\deg x :=e$ and 
 $\deg y:=e'$. Let $\fm=(x,y)$. Then $(S,\G,\fm)$ is a homogeneous coordinate ring for $\FF[p,q]$.

 The group homomorphisms 
 \begin{equation}
 \label{eq.actions}
\UseComputerModernTips
  \xymatrix{
 \CC^\times \ar[rr]^{\xi \mapsto (\xi^q,\xi^p)} &&  \G \ar[rrr]^{(\l_1,\l_2) \mapsto \l_1^p= \l_2^q} &&&  \CC^\times
 }
 \end{equation}
 induce homomorphisms
  \begin{equation}
 \label{eq.gr.gps}
 \UseComputerModernTips
  \xymatrix{
 \ZZ && \ar[ll]_{\phi'}  \G  && \ar[ll]_{\psi'}   \ZZ
 }
 \end{equation}
 between the rational character groups where $\phi'(ae+be')=aq+bp$ and $\psi'(1)=pe=qe'$.
 The group homomorphisms in (\ref{eq.actions}) induce morphisms
 \begin{equation}
 \label{morphisms}
\UseComputerModernTips
  \xymatrix{
\PP[q,p] \ar[rr]^\phiol  &&  \FF[p,q] \ar[rr]^\psiol && \PP^1
 }
 \end{equation}
 between the corresponding stack-theoretic quotients. The morphisms $\psiol$ and $\psiol\phiol$ are the natural morphisms to the coarse moduli spaces, and $\phiol$ is an isomorphism if and only if $(p,q)=1$
 (if and only if $\G$ is torsion-free). 
 
 The homomorphisms in (\ref{eq.gr.gps}) can also be interpreted as the natural maps  
  $$
\UseComputerModernTips
  \xymatrix{
\Pic \bigl(\PP[q,p] \bigr)   && \ar[ll]_{\phiol^*}  \Pic \bigl( \FF[p,q]  \bigr)&&\Pic \bigl( \PP^1\bigr) \ar[ll]_{\psiol^*}
 }
 $$
 between the Picard groups.  Similarly, the  homomorphisms in (\ref{eq.gr.gps}) induce the
 natural maps between the Grothendieck groups
  $$
\UseComputerModernTips
  \xymatrix{
K_0 \bigl(\PP[q,p] \bigr) \ar@{=}[d]   && \ar[ll]_{\phiol^*}  K_0 \bigl( \FF[p,q]  \bigr) \ar@{=}[d] &&  K_0 \bigl( \PP^1\bigr) \ar[ll]_{\psiol^*} \ar@{=}[d] 
\\
 {{\ZZ[u^{\pm 1}]}\over{((1-u^p)(1-u^q))}}  &&  {{\ZZ[s^{\pm 1},t^{\pm 1}]}\over{((1-s)(1-t),s^q-t^p)}} && {{\ZZ[v^{\pm p}]}\over{(v^p-1)^2}} 
 \\
 (u^p,u^q) && (s,t), \ar@{|->}[ll]  \quad  t^p=s^q && v \ar@{|->}[ll]
 }
 $$

 The morphisms in (\ref{morphisms}), and their associated inverse and direct image functors, are induced 
 by morphisms 
  $$
 \UseComputerModernTips
  \xymatrix{
(\CC[x,y],\ZZ,\fm)   &&  (\CC[x,y],\G,\fm) \ar[ll]_{(\phi,\phi')} &&    (\CC[x^p,y^q],\ZZ, (x^p,y^q))  \ar[ll]_{(\psi,\psi')}
 }
 $$
 between their homogeneous coordinate rings  where 
 \begin{itemize}
  \item 
  $\phi=\id_{\CC[x,y]}$ and $\psi$ is the inclusion $\CC[x^p,y^q] \to \CC[x,y]$, 
  \item 
  $\phi'$ and $\psi'$ are the homomorphisms between the grading groups in (\ref{eq.gr.gps}),
  \item 
  the $\ZZ$-grading on $\CC[x^p,y^q]$ is given by setting $\deg x^p=\deg y^q=1$ and 
  \item
   the $\ZZ$-grading on $\CC[x,y]$ is given by setting $\deg x= q$ and $\deg y=p$.   
\end{itemize}

  We call the point where $x$ vanishes {\sf north pole}.
The closed substack there is isomorphic to $B\mu_p$. We write $\cO_\fn$ for
the skyscraper sheaf at the north pole with the trivial equivariant structure. It fits
into an exact sequence 
 \begin{equation}
 \label{p.pres}
 \begin{CD}
 0 @>>>  \cO_{\FF[p,q]}(-e) @>{x}>> \cO_{\FF[p,q]} @>>> \cO_{\fn} @>>> 0
 \end{CD}
 \end{equation}
 and in terms of graded modules
 \begin{equation}
 \label{eq.On}
 \cO_\fn= \pi^* \biggl( {{\CC[x,y]}\over{(x)}}\biggr).
 \end{equation}
By (\ref{map.K0}), (\ref{loc.seq}), and (\ref{eq.On}), the class $[\cO_\fn]$   in $K_0(\FF[p,q])$ is  
$$
[\cO_\fn]= H_{{\CC[x,y]/(x)}} \cdot H_{\CC[x,y]}^{-1} = 1-t.
$$ 
 The other equivariant structures on $\cO_\fn$  are given by the modules 
 $$
 \cO_\fn(ie)=\pi^*\biggl({{\CC[x,y]}\over{(x)}}\biggr)(ie), \qquad i=0,1,\ldots,p-1,
$$
and $[\cO_\fn(-ie)]=t^i(1-t)$. The class in $K_0$ of a non-stacky point, say 
$\cO_\l=\pi^*\bigl(\CC[x,y]/(\l x^p-y^q)\bigr)$, $\l \ne 0$, is $1-t^p=1-s^q$.  We note that  
$$
[\cO_\l]=\sum_{i=0}^{p-1}  [\cO_\fn(-ie)].
$$
Because $1-t^p=t(1-t^p)=t(1-t^p)$ in $K_0$, 
twisting the structure sheaf of a non-stacky point does not change it.

{}

\end{document}